\documentclass[11pt]{amsart}
\usepackage{amsmath,amssymb}
\newcommand     {\comment}[1]   {}
\newcommand{\mute}[2] {}
\newcommand     {\printname}[1] {}
\newcommand{\labell}[1] {\label{#1}\printname{#1}}

\newcommand\fg {{\mathfrak g}}
\newcommand\fh {{\mathfrak h}}

\newcommand\ft {{\mathfrak t}}

\def    \oursetminus  {{\smallsetminus}}

\def    \fh     {{\mathfrak h}}

\def    \hG             {{\widehat{G}}}
\def    \hell   {\widehat{\ell}}

\def    \bs     {\boldsymbol}

\def    \to     {\longrightarrow}

\def    \tG     {\widetilde{G}}
\def    \tell   {\Tilde{\ell}}

\renewcommand{\Tilde}{\widetilde}
\usepackage{amssymb,latexsym}
\newcommand{\1}{{{\mathchoice {\rm 1\mskip-4mu l} {\rm 1\mskip-4mu l}
{\rm 1\mskip-4.5mu l} {\rm 1\mskip-5mu l}}}}

\newcommand{\Fmax}{F_{\rm max}}
\newcommand{\Fmin}{F_{\rm min}}

\newcommand{\R}{{\mathbb R}}

\newcommand{\Flux}{{\rm Flux}}
\newcommand{\Symp}{{\rm Symp}}

\newcommand{\id}{{\rm id}}

\newcommand{\GL}{{\rm GL}}

\newcommand{\Hh}{{\mathcal H}}

\newcommand{\Z}{{\mathbb Z}}
\newcommand{\C}{{\mathbb C}}

\newcommand{\om}{{\omega}}

\newcommand{\eps}{{\varepsilon}}
\newcommand{\de}{{\delta}}

\newcommand{\la}{{\lambda}}
\newcommand{\si}{{\sigma}}

\newcommand{\La}{{\Lambda}}

\newcommand{\Ham}{{\rm Ham}}

\newcommand{\MS}{{\medskip}}

\newcommand{\NI}{{\noindent}}

\newcommand{\QED}{\hfill$\Box$\medskip}

\newcommand{\CP}{{\mathbb{CP }}}

\newtheorem{theorem}{Theorem}[section]

\newtheorem{corollary}[theorem]{Corollary}

\newtheorem{thm}[theorem]{Theorem}

\newtheorem{example}[theorem]{Example}

\newtheorem{remark}[theorem]{Remark}

\newtheorem{lemma}[theorem]{Lemma}

\newtheorem{proposition}[theorem]{Proposition}

\begin{document}

\title{On nearly semifree circle actions}
\author{Dusa McDuff}
\thanks{First author partially supported by NSF grant DMS 0305939, and second by NSF grant DMS 0204448.}
\address{Department of Mathematics,
 Stony Brook University, Stony Brook, 
NY 11794-3651, USA}
\email{dusa@math.sunysb.edu}
\author{Susan Tolman}
\address{Department of Mathematics,
 University of Illinois at Urbana--Champaign, 
IL, USA}
\email{tolman@math.uiuc.edu}
\date{8 March 2005}

\maketitle

\MS

\begin{abstract}
Recall that an  effective circle action is semifree if the stabilizer 
subgroup of each point is connected.
We
show that if $(M, \om)$ is a coadjoint orbit of a compact Lie group $G$
then every element of
$\pi_1(G)$ may be represented by a semifree $S^1$-action.  
A theorem of McDuff--Slimowitz then implies that $\pi_1(G)$ injects into $\pi_1(\Ham(M, \om))$, which answers a 
question raised by Weinstein.   
We also show that a circle action on a manifold $M$ 
which is semifree near a fixed
point $x$ cannot contract in a compact
Lie subgroup $G$ of the diffeomorphism group unless the action is reversed 
by an element of $G$ that fixes the point $x$. 
Similarly, if a circle acts in a Hamiltonian fashion on a manifold
$(M,\omega)$ and the stabilizer of every point has at most
two components, then the  circle cannot contract in a compact Lie subgroup
of the group of Hamiltonian symplectomorphism unless the
circle  is reversed by
an element of $G$ 

\end{abstract}


\section{Introduction}

This paper is an attempt to understand topological properties of
Lie group actions.  Its starting point was the following theorem 
 of McDuff--Slimowitz~\cite{MSlim} concerning circle subgroups of 
$\Symp(M,\omega)$, 
the group of  symplectomorphisms of a symplectic manifold $(M,\om)$.
Recall that an effective  circle action 
is {\bf  semifree} if
the stabilizer subgroup of each point  in $M$ is either
the circle itself or the trivial group.
Also, we say that a circle subgroup $\La$ of a topological group $\Hh$ is {\bf 
essential in} ${\boldsymbol \Hh}$ if it represents a nonzero element in $\pi_1(\Hh)$ and 
{\bf inessential in} ${\bs \Hh}$ otherwise.

\begin{theorem}\labell{thm:sfr}  
Any semifree circle action on a closed symplectic manifold $(M,\omega)$
is essential in $\Symp(M,\om)$.
\end{theorem}

This is obvious if the
action is not Hamiltonian since in this case the flux homomorphism
$$
\Flux: \pi_1(\Symp(M,\om)) \to H^1(M,\R)
$$
does not vanish on $\La$.  However, if the action is Hamiltonian with
generating Hamiltonian $K:M\to \R$ then the result
is not so easy: the proof in~\cite{MSlim} involved studying
the Hofer length of the corresponding paths $\phi_t^K, t\in [0,T],$
in $\Ham(M,\om)$. 

The first result in  this paper uses the theorem above to answer
a question posed by Alan Weinstein in~\cite{Wei}.  
Let $G$ be a semisimple Lie group with Lie algebra $\fg$.
Let $M  \subset \fg^*$ be a coadjoint orbit,
together with the Kostant--Kirillov symplectic form $\omega$.
If the coadjoint action of $G$ on $M$ is effective, 
then $G$ is naturally a subgroup of $\Ham(M,\omega)$,
the group of Hamiltonian symplectomorphisms of $(M,\omega)$. 
This inclusion induces a natural map
from the fundamental group of $G$ to the fundamental group
of $\Ham(M,\omega)$. 
Weinstein
asks when this map is injective. 
We prove that this map is injective
for all compact semisimple Lie  groups.  
In~\cite{Vina} Vina established a 
special case of this result by quite different methods.

\begin{theorem}\label{cor:coadj}  
Let a compact semisimple Lie group $G$
act effectively on a coadjoint orbit $(M,\omega)$. 
Then the inclusion $G\to \Ham(M, \om)$ induces an 
injection from $\pi_1(G)$ to $\pi_1(\Ham(M,\om)).$
\end{theorem}

In view of  Theorem~\ref{thm:sfr}, this is an 
immediate consequence of the following result,
which we prove in Section \ref{ss:coadj}.

\begin{proposition}\label{WQ}
Let a compact semisimple Lie group $G$
act effectively on a coadjoint orbit $(M,\omega)$. 
Then every nontrivial element in $\pi_1(G)$ may be represented
by a circle that acts  semifreely on $M$. 
\end{proposition}

Theorem~\ref{thm:sfr} immediately implies that if a compact Lie group
$G$ acts effectively on a closed symplectic manifold $(M,\omega)$,
then any semifree circle
subgroup $\Lambda \subset G$ is essential in $G$.
The other results in the paper generalize this claim.

Observe that  Theorem~\ref{thm:sfr}
does not immediately extend  to 
the smooth (non-symplectic) category. For example, Claude
LeBrun pointed out to us that the circle action on $S^4$
induced by the diagonal
action of $S^1$ on $\C^2 = \R^4 \subset \R^5$ is semifree but gives a
nullhomotopic loop since $\pi_1(SO(5)) = \Z/2\Z$.
Nevertheless, the semifree condition does have consequences 
in the smooth category, even
if the action is  only semifree on
a neighborhood of a  component of the fixed point set; 
we shall say that such components are  
{\bf semifree}.  Further, given a circle subgroup
$\Lambda \subset G$ we say that $g \in G$ {\bf reverses}  
${\boldsymbol \Lambda}$ {\bf in} ${\boldsymbol G}$
if $g t g^{-1} = t^{-1}$ for all $t\in \La$.
Finally, a component $F$ of the fixed point set $M^\La$ of $\Lambda$ 
is {\bf symmetric in} ${\bs G}$
if there is an element $g \in G$ whose action on $M$ fixes $F$
pointwise and which reverses $\Lambda$.

\begin{thm}\label{thm:main0}
Let $\Lambda$ be a circle subgroup
of a  compact Lie group $G$ which acts effectively on a connected manifold $M$.
If there is a semifree component of the fixed point set $M^\Lambda$
which is not symmetric in $G$, then $\La$ is essential in $G$.
\end{thm}

\begin{example}\label{ex}\rm 
First, let $G = SU(2)$ act on $\CP^2$ by the defining representation on
the first two copies of $\C$,
and let $\Lambda \subset G$ be the circle subgroup
given by $\lambda \cdot [z_0:z_1: z_2]\mapsto
[\lambda z_0: \lambda^{-1} z_1: z_2]$.
This action has a  semifree fixed point, namely $[0,0,1]$.
Moreover,  this circle subgroup
is inessential in $G$.   
Therefore,  by the theorem above, there exists $g \in G$
which reverses the circle action and
fixes $[0,0,1]$.  In fact, we can take any $g$ which
lies in the normalizer  $N(\Lambda)$ but not in $\Lambda$ itself.
Note that $g^2 = -\1$ for any such $g$.

In contrast, consider the natural action
of $G= PU(3)$ on $\CP^2$, and
let $\Lambda \subset G$ be the circle subgroup
given by $\lambda \cdot [z_0:z_1: z_2]\mapsto
[\lambda^2 z_0: z_1: z_2]$.
This action is semifree and essential, 
but is not reversed by any $g \in G$.
To see this, note that the circle has order $3$ in $\pi_1(G)$,
whereas every circle that can be reversed has order $1$ or $2$.
\end{example}

One can weaken the semifree hypothesis in the above theorem, 
at the cost of adding
a global
isotropy assumption and working once more in the symplectic category.
We  say that a circle action
has {\bf at most twofold isotropy} if every point which is not  either fixed 
or free has stabilizer $\Z/(2).$  
Recall, also, that a symplectic action of $G$ on $(M,\om)$ is 
Hamiltonian if it is given by an equivariant moment map  $\Phi:M\to \fg^*$.

\begin{thm} \labell{thmtwo}
Let $\Lambda$ be a circle subgroup
of a  compact Lie group $G$ which acts effectively on a connected 
symplectic manifold $(M,\omega)$.
If $\Lambda$  has
at most twofold isotropy and if there is no  $g \in G$  
which reverses $\Lambda$, then
$\La$ is essential in $G$.
\end{thm}

\begin{example}\rm
This theorem does not extend to circle actions 
which have at most threefold isotropy.
For example, the action of $S^1$ on $\CP^3$ 
given by $\lambda \cdot [x,y,z,w] = [\lambda^2 x,  \lambda^{-1}y, 
\lambda^{-1}z,w]$ is inessential  in $PU(4)$.
However, since $\Fmax$ and $\Fmin$ are not diffeomorphic, 
this action has no reversor.

We also need the symplectic hypothesis.
To see this,  consider the obvious action of $SU(3)$ on $S^6: = 
\C^3\cup\{\infty\}$.  The subgroup $\La: = {\rm diag\,}(\lambda^2,
\lambda^{-1},\lambda^{-1})$ acts with at most twofold isotropy but 
has no reversor.  
\end{example}

\begin{remark}  \rm
If  $G$ is a simple group, we do not
need to assume that $(M,\omega)$ is symplectic in Theorem~\ref{thmtwo};
we only need to assume that there exists a point $p$ which
is fixed by a maximal torus containing $\Lambda$ but 
is not fixed by all of $G$.
Note that, in contrast, in the example above, 
the only points fixed by $\La$ are 
fixed by all of $G$.
\end{remark}

\begin{remark}\rm 
If $\Lambda$ is {\em any} circle subgroup of 
 $SO(3)$ -- or indeed a subgroup of any simple group
of type $B_n$, $C_n$, or $F_4$ --
then there exists a $g \in G$  which reverses $\Lambda$.   In this case,
Theorem~\ref{thmtwo} is trivial and the force of 
Theorem~\ref{thm:main0}
is that we can choose $g$ so that it also fixes $p$.
\end{remark}

\begin{remark}\rm
 In the proof of the above theorems, 
we pick a maximal torus $T$ which contains $\Lambda$.
The reversor $g$ that we construct lies in the normalizer $N(T)$ and
has the property that $g^2$ lies in  $T$. However, 
as we saw in Example~\ref{ex}, $g^2$ may not be equal to the identity. 
\end{remark}

Theorems \ref{thm:main0} and \ref{thmtwo}
have  the following easy corollaries:

\begin{corollary}
Consider a Hamiltonian
circle action $\Lambda$ on a closed symplectic manifold 
$(M,\omega)$ with moment map $K: M \to \R$,
normalized so that $\int_M K \omega^n = 0$.
If  $F$ is a semifree fixed component, then $\Lambda$
is essential in every compact subgroup
$G \subset \Symp(M,\omega)$ that contains it, unless  there is
a symplectomorphism $g$ of $M$ that fixes $F$ and reverses $\La$.
In this case, all the following hold:
\begin{enumerate}
\item $K(g(p)) = - K(p)$ for all $p \in M^{\La}$.
\item There is a one-to-one correspondence  between the
positive weights at $p$ and the negative weights at $g(p)$, and vice versa. 
\item $g$ induces an isomorphism on the image of the restriction
map in equivariant cohomology $H^*_{S^1}(M) \to H^*_{S^1}(M^\La)$.
\item $g(F) = F$.
In particular,
\begin{enumerate}
\item $K(F) = 0$.
\item The sum of the weights at $F$ is zero.
\end{enumerate}
\end{enumerate}
\end{corollary}

\begin{corollary}
Consider a Hamiltonian
circle action $\Lambda$ on a closed symplectic manifold 
$(M,\omega)$ with moment map $K: M \to \R$,
normalized so that $\int_M K \omega^n = 0$.
If the action has at most twofold isotropy,  then $\Lambda$
is essential in every compact subgroup
$G \subset \Symp(M,\omega)$ that contains it, unless  there is
a symplectomorphism $g$ of $M$ that reverses $\La$.  In this case,
 all the following hold:
\begin{enumerate}
\item $K(g(p)) = - K(p)$ for all $p \in M^{\La}$.
\item There is a one-to-one correspondence  between the
positive weights at $p$ and the negative weights at $g(p)$, and vice versa. 
\item $g$ induces an isomorphism on the image of the restriction
map in equivariant cohomology $H^*_{S^1}(M) \to H^*_{S^1}(M^\La)$.
\end{enumerate}
\end{corollary}

It is unknown whether the existence of such $g$ is necessary for $\La$ to be inessential in $\Symp(M,\om)$. 
We make partial progress towards answering this question in
\cite{MT}.

All the results in this paper are proved by 
a case by case study of the structure of semisimple Lie algebras.

\section{Coadjoint orbits}\label{ss:coadj}

In this section, we prove Proposition~\ref{WQ}.
We begin with a brief review of a few facts about Lie groups.

Each simply connected compact semisimple Lie group is a product of simple 
factors, and its center is the product of the centers of its
simple factors.   Moreover,
since its Lie algebra splits into a corresponding sum, the 
coadjoint orbits also are products of coadjoint orbits of simple groups. 
Therefore, we may assume that $G$ is simple. 

Let $G$ be a compact simple Lie group.
Let $\tG$ denote the universal cover of $G$, and
$\hG$ denote the quotient of $G$ by its center.
Let $\fg$ denote the Lie algebra of $G$, and
let $\ft$ denote  the Lie algebra of a maximal torus
$T \subset G$. 
Let $\ell \subset \ft$, $\tell \subset \ft$, and $\hell \subset \ft$ be the lattices consisting of
vectors $\xi \in \ft$ whose exponential is the identity
in $G$, $\tG$, and $\hG$, respectively.
There is a one-to-one correspondence between $\ell$ and circle subgroup
of $G$, $\tell$ and circle subgroups
of $\tG$, and $\hell$ and circle subgroups of $\hG$,
given by sending $\lambda$ to $t \rightarrow \exp(t\lambda)$.
Note that $\tell \subseteq \ell \subseteq \hell$. 
Because $\tG$ is simply 
connected,
$$
\pi_1(G) \;\cong \;\ell/\tell \;\subseteq \;\hell/\tell \;\cong \; \pi_1(\hG).
$$

Let $\ft^*$ denote the dual to $\ft$, and
let $\Delta \subset \ft^*$ denote the set of {\bf roots} of $G$, i.e. the 
nonzero weights of the adjoint action  $T$ on $\fg_\C$,
where $\fg_\C$ is the complexification of $\fg$.
The lattice $\hell$ is dual to the lattice in $\ft^*$ generated by the 
roots, i.e. $\la\in \hell$ precisely when $\eta(\la)\in \Z$ for all 
$\eta\in \Delta$.
If we use the Killing form $(\cdot, \cdot)$  to identify
$\ft$ and $\ft^*$, then  $\tell$ is generated by the set
$$
\left\{ \left. \frac{2 \eta}{(\eta,\eta)}\ \right| \eta \in \Delta \right\}.
$$
Further the set of weights at any fixed point  $p$ 
for the action of $T$ on $M$ is a nonempty subset of the set of roots.
Therefore
the result will follow if we find a representative $\lambda$
for each nontrivial class in $\hell/\tell$ such that 
$|\eta(\lambda)| \leq 1$ for every $\eta \in \Delta$.

We will check this on a case by case basis;
in each case we will use the Killing
form to identify $\ft$ and $\ft^*$.
Let $( \cdot, \cdot)$ be the standard metric on $\R^k$
with the standard basis $e_1,\dots,e_k$, and define
$$
\epsilon_i = e_i - \frac{1}{k} \sum_{j=1}^k e_j.
$$
\MS

\NI{\bf (I)}\,\,
For the group $A_n$, where $n \geq 1$, 
$\ft = \ft^* = \left\{ \lambda \in \R^{n+1} \left| \ \sum \lambda_i = 0 \right. \right\}$
and the roots are $\epsilon_i - \epsilon_j = e_i - e_j$ for $i \neq j$.  
Hence
$$
\hell = \{ \lambda \in \ft \mid \lambda_i - \lambda_j \in \Z \ \forall \ i,j \},
\ \mbox{and} \quad 
\tell =   \ft \cap \Z^{n+1}.
$$
As representatives for the quotient $\hell/\tell \cong \Z/(n+1)$, we  take
$\lambda =   \sum_{i=1}^k \epsilon_i$ for $0 \leq k \leq n$.
\MS

\NI{\bf (II)}\,\,
For the group $B_n$, where $n \geq 2$, $\ft^* = \R^n$ and
the roots are $\pm e_i$ and $\pm e_i \pm  e_j$ for $i \neq j$.
Hence
$$
\hell = \Z^n,\ \mbox{and} \qquad
\tell = \left\{ \lambda \in \Z^n 
\left| \  \sum \lambda_i \in 2\Z \right. \right\}.
$$
As representatives of the quotient $\hell/\tell \cong \Z/(2) $,
we  take $0$ and $e_1$.
\MS

\NI{\bf (III)}\,\,
For the group $C_n$, where $n \geq 3$,  $\ft^* = \R^n$ 
and the  roots are  $\pm 2 e_i$ and $\pm e_i \pm  e_j$ for $i \neq j$.
Hence
$$
\hell = \{ \lambda \in  \R^n \mid  \lambda_i \pm \lambda_j \in \Z,
\ \forall \ i,j  \},
\ \mbox{and} \qquad
\tell = \Z^n.
$$
As representatives of the quotient $\hell/\tell \cong \Z/(2)$ ,
we  take $0$ and $\frac{1}{2} \sum_{i=1}^n e_i$.
\MS

\NI{\bf (IV)}\,\,
For the group $D_n$, where $n \geq 4$, $\ft^* = \R^n$ and the 
roots are  $\pm e_i \pm e_j$ for $i \neq j$.  Hence
$$
\hell =
\{ \lambda \in  \R^n \mid  \lambda_i \pm \lambda_j \in \Z, \ \forall \ i,j  \},
\ \mbox{and} \qquad
\tell = \left\{ \lambda \in \Z^n 
\left| \  \sum \lambda_i \in 2\Z \right. \right\}.
$$
The quotient $\hell/\tell$ is isomorphic to
$\Z/(2) \oplus \Z/(2)$ if $n$ is even, and
to $\Z/(4)$ if $n$ is odd.
Either way, as representatives of $\hell/\tell$, we  take
$0$, $e_1$, $\frac{1}{2}\sum_{i=1}^n e_i$
and $\frac{1}{2}\sum_{i=1}^n e_i - e_n$.
\MS

\NI{\bf (V, a)}\,\,
For the group $E_6$,  $\ft^* = \R^6$
and the roots  are $2 \epsilon$,  $\epsilon_i - \epsilon_j$,
and $\epsilon_i + \epsilon_j + \epsilon_k \pm \epsilon$
for $i, j,$ and $k$ distinct,
where $\epsilon = \frac{1}{2 \sqrt{3}}(1,1,1,1,1,1)$.
Hence, 
$$
\hell =
\Big\{ n \epsilon + (\xi_1,\ldots,\xi_6) \in  \R^6  \Big| 
\sum_{i=1}^6 \xi_i = 0, n \in \Z, \
\frac{n}{2} + 3 \xi_i \in  \Z \ \mbox{and} \
\xi_i - \xi_j \in \Z \ \forall \ i,j  \Big\},
\ \mbox{and} $$
$$
\tell =
\Big\{ n \epsilon + (\xi_1,\ldots,\xi_6) \in  \R^6   \Big|
\sum_{i=1}^6 \xi_i = 0, n \in \Z \ \mbox{and} \ 
\frac{n}{2} +  \xi_i \in  \Z \ \forall \ i  \Big\}.
$$
As representatives of the quotient $\hell/\tell \cong \Z/(3)$ ,
we  take $0$, $\epsilon_1 + \epsilon_2$, and $-\epsilon_1 - \epsilon_2$.
\MS

\NI{\bf (V, b)}\,\,
For the group $E_7$, 
$\ft = \ft^* = \left\{ \lambda \in \R^{8} \mid \sum \lambda_i = 0 \right\}$,
and the roots are $\epsilon_i - \epsilon_j$, and
$\epsilon_i + \epsilon_j + \epsilon_k + \epsilon_l$ for
$i, j, k,$ and $l$ distinct.  Hence 
$$\hell =
\{ \lambda \in \ft \mid
4 \lambda_i \in \Z \ \mbox{and} \  \lambda_i - \lambda_j \in \Z\  \forall\  i, j \},
\ \mbox{and} \quad \tell  
= \{ \lambda \in \ft \mid
 \lambda_i \pm \lambda_j \in \Z\ \forall \ i, j
\}.
$$
As representatives for the quotient $\hell/\tell \cong \Z/2\Z$, we  take
$0$ and $\epsilon_1 + \epsilon_2$.\MS

Every group of type  $E_8$, $F_4$, and $G_2$ is simply connected,
so no further argument is necessary.
\QED

\section{Lie Group Actions}\label{sec:lie}

This section contains proofs of Theorems \ref{thm:main0} and \ref{thmtwo}. We begin by stating a lemma about root systems, that is proved at the end.  We shall 
 always assume that the  positive Weyl chamber is closed.

\begin{lemma}\labell{claims}
Let $G$ be a simply connected compact simple Lie group.
Let $\ft$ be the Lie algebra  of a maximal torus $T \subset G$.
Let $\tell$ be the integral lattice,  let $\Delta$ denote the set of roots,
and let $W$ denote the Weyl group. Use the Killing form to identify $\ft$ and
$\ft^*$.
Fix $\lambda \in \tell$.
Choose a positive Weyl chamber which contains $\lambda$.
Let $\delta \in \Delta$ denote the highest root.
Then the following claims hold:
\smallskip

\NI{\rm (a)}
 If $(\lambda,\delta) \leq 2$, 
then there exist orthogonal roots $\eta_1,\ldots,\eta_k \in \Delta$
so that $\lambda = \sum a_i \eta_i$ and so that
$(\lambda,\eta_i) = a_i (\eta_i,\eta_i) = 2$ for all $i$.
\smallskip

\NI{\rm (b)}
Let $L \subset \Delta$ be a set of roots 
which contains every root $\eta \in \Delta$ 
such that $\delta + \eta$ or $\delta - \eta$ is also a root.  Assume
also that $L$ is  closed under addition,  that is, it  
contains every root which can be written as the sum of
roots in $L$.  Then $L$ contains all roots.
\smallskip

\NI{\rm (c)}
 If $(\lambda,\delta) > 2$ and $-\id: \ft \to \ft$ is not an 
element of the Weyl group,
then for every nonzero weight $\alpha \in \tell^*$
there  exists $\sigma \in W$ so that $|(\sigma \cdot \alpha,\lambda)| > 1$.
\smallskip

\NI{\rm (d)}
 If $-\id: \ft \to \ft$ is not an 
element of the Weyl group,
then $\delta$ is the only root which lies in the positive Weyl chamber.
\end{lemma}

Using this result, we can find find elements which reverse certain
circle subgroups of simply connected compact simple Lie groups.
Note that because $G$ 
is simply connected,  every circle subgroup of $G$ is inessential in $G$.

\begin{lemma} \labell{le:simple}
Let $\Lambda $ be a circle subgroup of a simply connected compact
simple Lie group $G$.
\smallskip

\NI{\rm (i)}
Let $\rho : G \to \GL(V)$ be a nontrivial representation of $G$.
If $\Lambda$  acts semifreely on $V$ then there exists
$g \in G$ that reverses $\Lambda$.
\smallskip

\NI{\rm (ii)}
Let $H \varsubsetneq G$ be a proper subgroup containing $\Lambda$.
If the adjoint action of $\Lambda$ on $\fg/\fh$ is semifree,
then there exists  $h \in H$  that reverses $\Lambda$.
\smallskip

\NI{\rm (iii)}
Let $H \varsubsetneq G$ be a proper subgroup containing a maximal
torus which contains $\Lambda$.
If the natural action of $\Lambda$ on $G/H$ has at most twofold isotropy,
then there exists $g \in G$ that reverses $\Lambda$.
\end{lemma}

The assumption in (ii) above is a special case of (i) since
the representation $V$ is restricted; however,
the conclusion is stronger since it asserts that the reversor lies 
in $H$.   Statement (ii) and (iii) are also related: the former makes a strong 
assumption about
the action induced by $\La$ on the tangent space to $G/H$ 
at the fixed point $eH$, the latter  makes a weaker assumption about the 
action at all the fixed points on $G/H$.

We will now use the claims in Lemma \ref{claims} to prove Lemma \ref{le:simple}.
Let $T$ be a maximal torus which contains $\Lambda$.
Let $\tell \subset \ft$ denote the integral  lattice.
Let $\la \in\tell$  be the vector corresponding to $\Lambda$.
Choose a positive Weyl chamber which contains $\lambda$.
Let $\delta \in \Delta$ denote the highest root.

Recall that the  Weyl group $W$ is the quotient $N(T)/T$, where $N(T)$ is the
the normalizer of $T$ in $G$.
Every root $\eta$ gives rise to an element $w_\eta \in W$
whose action on $\ft^*$ is given by
$w_\eta(\beta) = \beta - \frac{2(\eta,\beta)}{(\eta,\eta)} \eta$.
\MS

\NI
{\bf Proof of Lemma~\ref{le:simple} (i).}
Let $\rho : G \to \GL(V)$ be a nontrivial representation of $G$.
Assume that  $\Lambda$  acts semifreely on $V$.

Suppose first that $(\lambda,\de) \leq 2$.
By claim (a), there exist
orthogonal roots $\eta_1,\ldots,\eta_k \in \Delta$ so that
$\lambda = \sum a_i \eta_i$.
Since the roots are orthogonal,
for each $\eta_i$ the associated element of the Weyl group $w_{\eta_i}$
takes $\eta_i$ to $-\eta_i$ and leaves $\eta_j$ fixed for
all $j \neq i$.  Hence, their product
$w = w_{\eta_1} \cdots w_{\eta_n}$ takes 
$\lambda$ to $-\lambda$, and so reverses $\La$.

So assume instead that 
$(\lambda,\de) > 2$.
If $-\id$ is in the Weyl group, then statement (i) is trivial.
So we assume that it is not.
Let $T$ act on $V$ via restriction,
and pick any nonzero weight $\alpha \in \tell^*$ in the weight
decomposition.  
By claim (c), we can find some $\sigma \in W$
such that $|(\sigma \cdot \alpha, \lambda) | >1.$
Since $\sigma \cdot \alpha$ also
appears in the weight decomposition,  
this contradicts the assumption that the action of $\Lambda$ on $V$
is semifree.\QED

\NI
{\bf Proof of Lemma~\ref{le:simple} (ii).}
Let $H \subseteq  G$ be a proper subgroup which contains $\Lambda$.
Assume that the adjoint action of $\Lambda$ on $\fg/\fh$ is semifree.
Let $L$ be the set of roots $\eta \in \Delta$ so that the
associated weight space $E_\eta \subset \fg_\C$ lies in $\fh_\C$.
Clearly, if $|(\eta, \lambda)| > 1$, then $\eta \in L$.

Suppose first that $(\lambda,\de) \leq 2$.
By claim (a), there exist orthogonal
roots $\eta_1,\ldots,\eta_k \in \Delta$
so that  $\lambda = \sum a_i \eta_i$ and so
that $(\lambda,\eta_i) = 2$ for every $i$.
Since  $(\eta_i,\lambda) = 2$,
$\eta_i$ lies in $L$ for all $i$.
Hence, the associated element of the Weyl group
$w_{\eta_i}$ lies in $H$ for all $i$.
Thus $w = w_{\eta_1} \cdots w_{\eta_k}$ must lie in $H$.

So assume instead that  $(\lambda,\delta) > 2$. 
We see immediately that  $\delta$ and $-\delta$ lie in $L$.
If $\eta$, $\eta'$ and $\eta + \eta'$
are all roots, then $[E_\eta,E_\eta'] = E_{\eta + \eta'}$.
Hence, since $\fh_\C$ is closed under Lie bracket,
if $\eta$ and $\eta'$ are in $L$ then $\eta + \eta' \in L$ also,
that is, $L$ is closed under addition.
Additionally, if $\eta$ and $\eta'$ are
roots such that $\delta = \eta  + \eta'$, then either
$(\lambda,\eta) > 1$ or $(\lambda,\eta') > 1$.
If the former holds, then $\eta$ and $-\eta$  lie in $L$.
Since $L$ is closed under addition, so do $\eta' $ and $-\eta'$.
The other case is identical.
Thus,  claim (b) implies that every root lies in $L$.
This contradicts the claim that $H$ is a proper subgroup.\QED

\NI
{\bf Proof of Lemma~\ref{le:simple} (iii).}
Let $H \subsetneq G$ be a proper subgroup which contains the 
maximal torus $T$, and assume that the
natural action of 
$\La\subset T$
on $G/H$ has at most twofold isotropy.

If $(\lambda,\delta) \leq 2$,  then part (iii) 
follows by the argument used to prove part (i). 
So assume that $(\lambda, \delta) > 2$.
We may also assume that
$-\id$ does not lie in the Weyl group,
because otherwise the claim is trivial.
Since $H \subset G$ is proper, there exists at least one root
$\eta$  so that the associated weight space $E_\eta$ is not contained in $\fh_\C$.
Then there is $\sigma \in W$ so that the root
$\sigma \cdot \eta$ lies in the 
positive Weyl chamber.  Hence by  (d) $\sigma \cdot \eta = \de$, and so 
$|(\sigma \cdot \eta,\lambda)| > 2.$
Choose $\tilde{\sigma} \in N(T)$
which descends to $\sigma$.  Then $\tilde{\sigma} H$ is 
a fixed point  for $T$, and $\sigma \cdot \eta$ is one of
the weights for $\La$
at this fixed point.
This contradicts the fact that the action has at most twofold isotropy.\QED

We are now ready to  deduce Theorems~\ref{thm:main0} and~\ref{thmtwo}.
In both cases, we will do this by proving the contrapositive,
that is, we will assume that $\Lambda$ is an inessential circle subgroup
and use this to construct a reversor.

Let $\tG$ denote the universal cover of $G$. Then
$\tG$ is the direct product of a compact simply connected 
semisimple Lie group and a vector space.
Since $\Lambda$ is inessential, it
lifts to a circle subgroup of  $\tG$.
Since this lift must lie in the compact part of $\tG$, we 
may assume without loss of generality that
$\tG$ is a compact simply connected semisimple Lie group.

In fact, it is enough to prove these claims  for the universal
cover of $G$, as long as we no longer insist on an effective action
but instead allow a finite number of elements of the group to act trivially
on $M$. 
Thus we may assume that $G$ is the product of compact simple and
simply connected groups
$G_1 \times \cdots \times G_n$. 
Let $\Lambda_i$ be the projection of $\Lambda$ to $G_i$.
Without loss of generality, we may assume that $\La_i\ne \{\id\}$ for 
all $i$.
\MS

\NI
{\bf Proof of Theorem~\ref{thm:main0}.}
Let $G = G_1\times \cdots\times G_n$ as above.
Choose $p\in F$ and let
$H \subset G$ be the stabilizer of $p$.  Then $\La\subset H$.  
There exists a representation $V$ of $H$, called the {\bf isotropy
representation}, so that a neighborhood of the $G$-orbit
through $p$ is equivariantly diffeomorphic to 
a neighborhood of the zero section of $G \times_H V$.
Fix some simple factor $G_i$, and let $H_i = H \cap G_i$.

Assume first that $H_i$ is a proper subgroup.
Note that $\fg_i$  is invariant under the action of $\Lambda$.
Thus, since $\Lambda$ acts semifreely on $\fg/\fh$ via the adjoint action,
$\Lambda_i$ acts semifreely on  $\fg_i/\fh_i$.
Thus, by  Lemma~\ref{le:simple} (ii)
there exists an element  $h_i  \in H_i$ that reverses $\Lambda_i$.

So assume on the contrary that $H_i = G_i$.
Let $\Lambda'$ be the projection of $\Lambda$ onto
the product of all the simple factors except $G_i$.
Since $\La_i\subset G_i\subset H$ and $\La\subset H$,
we must have $\La'\subset H$.  Hence $\La'$ acts on $V$.
For any integer $k$, let $V_k$ denote the subspace of
$V$ on which $\Lambda'$ acts  with weight $k$.
Since $\Lambda'$ commutes with $G_i$,  $V_k$
is a representation of $G_i$.
Since only a finite number of elements of $G$ act trivially
on $M$, $G_i$ must act nontrivially on
$G \times_H V$, and hence also on $V$.
Therefore,  there is some $k$ so that the 
representation of $G_i$ on $V_k$ is nontrivial.
Because $G_i$ is simple,
$\Lambda_i$ must act with both  positive
and negative weights on $V_k$.
But the weights for the action of $\Lambda$ on
$V_k$ are the weights for the action of $\Lambda_i$ shifted by $k$.
Hence,
because $F$ is a semifree fixed point component,
$k = 0$ and the action of $\Lambda_i$
on $V_k$ is itself semifree. 
Therefore by  Lemma~\ref{le:simple} (i)  there exists  
$h_i \in G_i = H_i$ that reverses $\Lambda_i$.

Since $h_i$ reverses $\Lambda_i$ for each
$i$, $g = (h_1,\ldots,h_n)$ reverses  $\Lambda$, as required.
Moreover, since
$H_1 \times \cdots \times H_n\subset H$
(in general they are not equal), $g$
lies in $H$, and hence fixes $p$.
\QED

\NI
{\bf Proof of Theorem~\ref{thmtwo}.}
Fix some simple factor $G_i$.
Let $W$ be the Weyl group of $G_i$.
Let $T \subset G_i$ be a maximal torus of $G_i$ containing $\Lambda_i$.
Let  $\Phi: M \to \ft^*$ be the moment map for the $T$-action.
Pick any $\xi \in \ft$ so that the one parameter subgroup
generated by $\xi$ is dense in $T$.
Let $p$ be  any point which maps to 
the minimum value of  $\Phi^\xi$, the
component of $\Phi$ in the direction $\xi$.
By construction, $p$ is a fixed point for $T$.
Assume first that  $\Phi^\xi(p) = 0$, that is,
the function $\Phi^\xi$ is nonnegative on
$M$. Since the moment polytope $\Phi(M)$ is 
invariant under the Weyl group $W$,
this implies that
$\Phi^{\si \cdot \xi}$ is  also nonnegative on $M$
for all $\si\in W$.
Because $G_i$ is simple and $\xi$ is a generic point of $\ft$,
for any nonzero $x\in \ft^*$ there exists an 
element $\si\in W$ such that $(\si\cdot\xi,x) < 0$.
Applying this to $x\in \Phi(M)\oursetminus \{0\}$, we see that
$\Phi(M)$ must be the single point $\{0\}$,
which is impossible, because the action is effective.
Therefore, $\Phi(p) \neq 0$. 

Now let us reconsider the action of $G$ on $M$.
Let $H$ be the stabilizer of $p$ in $G$, and let $H_i = H \cap G_i$.
Since $\Lambda$ acts with 
at most twofold isotropy on 
$G/H\subset M$,  $\Lambda_i$
acts with at most twofold  isotropy on $G_i/H_i$.
Since $\Phi(p)$ is not zero, the stabilizer of $\Phi(p)$ in $G_i$.
is a proper subgroup of $G_i$.  Since $\Phi$ is equivariant,
this implies that $H_i$ is a proper subgroup of $G_i$.
By Lemma \ref{le:simple} (iii), this implies that
there exists $g_i \in G_i$ which reverses $\Lambda_i$.
Then $(g_1,\ldots,g_n)$ reverses $\La$.
\QED

\NI
{\bf Proof of Lemma \ref{claims}.}
We now prove claims (a)-(d) on a case by case basis,
using the classification of compact simple Lie groups.
We will use the notation of 
\S\ref{ss:coadj}.
Note, however, that here $G =\tG $ since $G$  is simply connected.
\MS

\NI{\bf (I)}\,\,
Recall that for the group $A_n$, where $n \geq 1$,  $\ft = \ft^*
= \{ \xi \in \R^n \mid \sum \xi_i = 0 \}$, the
roots are $\epsilon_i - \epsilon_j
= e_i-e_j$ 
for $i \neq j$, and 
the integral lattice is $\tell = \Z^{n+1} \cap \ft$. 
The positive Weyl chamber is\footnote
{For uniformity, we shall always
use the lexigraphical order to choose the positive Weyl chamber.}
 $$
 \{\xi \in \ft \mid
\xi_1 \geq \cdots \geq \xi_{n+1}\}.
$$
The highest root is $\delta = e_1-e_{n+1}$.

If $(\lambda,\delta) = \lambda_1 - \lambda_{n+1} \leq 2$, 
then $|\lambda_i| \leq 1$
for all $i$. Since $\sum_i \lambda_i = 0$ and  $\lambda_i \in \Z$ for all $i$,
there are an equal number of $+1$'s and $-1$'s, and the rest are $0$'s.
Hence, $\lambda$ is the sum of orthogonal roots of the form
$\eta = e_i - e_j$.  Since $(\eta,\eta) = 2$,
this proves claim (a).

Since $\delta = (e_1 - e_k) + (e_k - e_{n+1})$,
the roots  $\pm (e_1 - e_k)$ and $\pm (e_k - e_{n+1})$ lie in $L$ for 
all $1 < k < n+1$.
If neither $i$ nor $j$ is equal  equal to $1$, then
$e_i - e_j = -(e_1 - e_i) +  (e_1 - e_j)$ is also in $L$.
This proves claim (b).

We now prove (c).
The weight lattice  is
$\tell^* = \{ \alpha \in \ft \mid \alpha_i - \alpha_j \in \Z \  \;
\forall \  i,  j \}$.
By permuting the coordinates of $\alpha$,
we may assume $\alpha_1 \geq \cdots \ge \alpha_n$.
Since $\alpha \neq 0$, there exists $k \in (1,\ldots,n)$
such that $\alpha_{k} - \alpha_{k+1} >0$; since
this difference lies in $\Z$, it must be at least $1$.
Since $\lambda_1 - \lambda_{n+1} = \lambda_1 + \sum_{i=1}^n \lambda_i > 2$ and 
$\lambda_i \geq \lambda_{i + n - k}$, 
$\sum_{i=1}^k \lambda_i + \sum_{i = 1}^{n+1 - k} \lambda_i > 2$.
Therefore, either
$\sum_{i=1}^k \lambda_i >1$ or $ \sum_{i =1}^{n+1-k} \lambda_i > 1$.
In the former case,
$$
(\alpha,\lambda) = \sum_{j=1}^n \left( (\alpha_j - \alpha_{j+1}) \sum_{i=1}^j
\lambda_i \right) \geq (\alpha_{k} - \alpha_{k+1}) \sum_{i=1}^k \lambda_i > 1.
$$
In the latter case, let $\alpha'$ be obtained from $\alpha$ by 
the permutation which 
reverses the coordinates, so that $\alpha'_i = \alpha_{n + 2 - i}$.
Then
$$
(\alpha',\lambda) = \sum_{j=1}^n \left( (\alpha_{n + 2 - j} - \alpha_{n + 1 - j})
\sum_{i=1}^j
\lambda_i \right) \leq (\alpha_{k + 1} - \alpha_{k}) \sum_{i=1}^{n + 1 - k} 
 \lambda_i < -1.
$$
The only facts we have used are that 
$\ft = \{ \xi \in \R^n \mid \sum_i\xi_i = 0 \}$,
that the Weyl group contains the permutation group $S_n$, and that 
$\alpha_i - \alpha_j \in \Z$ for any $\alpha \in \tell^*$. 

Finally,  $\delta$ is the only root in the positive Weyl chamber.

\MS

\NI{\bf (II)}\,\,
Recall that for the group $B_n$, where $n \geq 2$, 
$\ft = \ft^* = \R^n$, the roots are $\pm  e_i$ and $\pm e_i \pm e_j$ for $i \neq j$,
and the integral lattice is $\tell = \{ \xi  \in \Z^n \mid
\sum_i \xi_i \in 2 \Z \}.$
The positive Weyl chamber is
$\{\xi \in \ft \mid \xi_1 \geq \cdots \geq \xi_n \geq 0 \}$.
The highest root is $\delta = e_1 + e_2$.

If $(\lambda,\delta) = 
\lambda_1 + \lambda_2 \leq 2$, then either
$\lambda_1 = 2$ and $\lambda_i = 0$ for all $i 
\neq 1$, or $\lambda_i \leq 1$ for all $i$.
Either way, since $\sum_i \lambda_i \in 2 \Z$, 
we can write $\lambda$ as the sum 
of orthogonal roots $\eta_i$ such that $(\eta_i,\eta_i) = 2$.

Since $\delta = (e_1 - e_k) + (e_2 + e_k) = (e_1 + e_k) + (e_2 - e_k)$,
the roots $\pm e_1 \pm  e_k$ and $\pm  e_2 \pm  e_k$
lie in $L$ for $k \neq 1 $ or $2$.
Since $\delta = (e_1) + (e_2)$,  the roots $\pm e_1$ and $\pm e_2$ lie in $L$.
Every root can be written as a sum of these roots.

Since $-\id$ lies in the Weyl group, we are done.
\MS

\NI{\bf (III)}\,\,
Recall that for the group  $C_n$, where  $n \geq 3$, $\ft = \ft^* = \R^n$,
the roots are $\pm 2 e_i$ and $\pm e_i \pm e_j$ for $i \neq j$, and the
integral lattice is $\tell = \Z^n$.
The positive Weyl chamber is
$\{\xi \in \ft \mid \xi_1 \geq \cdots \geq \xi_n \geq 0\ \}$.
The highest root is $\delta = 2 e_1$.

If $(\lambda,\delta) = 2 \lambda_1 \leq 2$, then
$\lambda_i \leq 1$ for all $i$.
Since $\lambda \in  \Z^n$, we can write $\lambda$
as half the  sum of orthogonal roots of the form $ 2 e_i$. 
Note that $(\lambda,   2 e_i) = 2$.

Since $\delta = (e_1 - e_k) + (e_1 + e_k)$,
the roots $\pm e_1 \pm  e_k$ 
lie in $L$ for $k \neq 1 $.
Every root can be written as a sum of these roots.

Since $-\id$ lies in the Weyl group, we are done.

\MS

\NI{\bf (IV)}\,\,
Recall that for the group $D_n$, where $n \geq 4$, 
$\ft = \ft^* = \R^n$, the roots are $\pm e_i \pm e_j$ for $i \neq j$,
and the integral lattice is
 $\tell =
\{ \xi \in  \Z^n \mid \sum \xi_i  \in 2\Z \}.$
The positive Weyl chamber is
$\{ \xi \in \ft \mid \xi_1 \geq \cdots \geq \xi_{n-1} \geq |\xi_n|  \}$.
The highest root is $\delta = e_1 + e_2$.

If $(\lambda,\delta) = 
\lambda_1 + \lambda_2 \leq 2$, then either
$\lambda_1 = 2$ and $\lambda_i = 0$ for all $i 
\neq 1$, or $|\lambda_i| \leq 1$ for all $i$.
Either way, since $\sum_i \lambda_i \in 2 \Z$, 
we can write $\lambda$ as the sum 
of orthogonal roots $\eta_i$ such that $(\eta_i,\eta_i) = 2$.

Since $\delta = (e_1 - e_k) + (e_2 + e_k) = (e_1 + e_k) + (e_2 - e_k)$,
the roots $\pm e_1 \pm  e_k$ and $\pm  e_2 \pm  e_k$
lie in $L$ for $k \neq 1 $ or $2$.
Every root can be written as a sum of these roots.

Now assume that $(\delta,\lambda) = \lambda_1 + \lambda_2 > 2$.
Consider a nonzero weight $\alpha \in \tell^* = \{ \alpha \in \R^n \mid
\alpha_i \pm \alpha_j \in \Z \ \forall\ 
i,  j \}.$
By applying the Weyl group, we may assume $\alpha$ lies in the positive Weyl chamber.
Since $\lambda$ also lies in the positive Weyl chamber,
$\alpha_i \lambda_i \geq 0$ for all $i \neq n$.
Moreover,  since $\alpha_{n-1} \geq |\alpha_n|$, and 
$\lambda_{n-1} \geq |\lambda_n|$, $\alpha_{n-1} \lambda_{n-1} +
\alpha_n \lambda_n \geq 0$.
Therefore, $\alpha_3 \lambda_3 + \cdots + \alpha_n \lambda_n \geq 0.$
(Here, we have used that $n \geq 4$.) 
Since $\alpha$ is nonzero, either $\alpha_1 \geq 1$, 
or $\alpha_1 = \alpha_2 = \frac{1}{2}$.
In either case, $\alpha_1 \lambda_1 + \alpha_2 \lambda_2 > 1$ .
(In the first case, we use the fact  that $\lambda_1 + \lambda_2 > 2$ and
$\lambda_1 \geq \lambda_2$ implies that $\lambda_1 > 1$.) 
Therefore,
$(\alpha,\lambda) \geq \alpha_1 \lambda_1 + \alpha_2 \lambda_2 > 1$,
This proves claim (c).

Finally,  $\delta$ is the only root in the positive Weyl chamber.

\MS

\NI{\bf (V, a)}\,\,
Recall that  for the group $E_6$, $\ft = \ft^* = \R^6$
and the roots  are $2 \epsilon$,  $\epsilon_i - \epsilon_j$,
and $\epsilon_i + \epsilon_j + \epsilon_k \pm \epsilon$
for $i, j,$ and $k$ distinct,
where $\epsilon = \frac{1}{2 \sqrt{3}}(1,1,1,1,1,1)$. Therefore
$$
\tell = \Big\{ n\epsilon + (\xi_1,\ldots,\xi_6)\; \Big|\; \sum_{i=1}^6 \xi_i = 0, n \in \Z,
\mbox{and} \ \frac{n}{2}+ \xi_i \in \Z \ \forall \ i \Big\}.
$$
The positive Weyl chamber is
$$
\Bigl\{ n \epsilon + (\xi_1,\ldots,\xi_6) \in \ft\; \Big|  \;
\sum_{i=1}^6 \xi_i = 0,  \xi_2 \geq \cdots \geq \xi_6,\; \xi_1 + 
\xi_5 + \xi_6 \geq n/2 \geq 0\Bigr\}.
$$
(Note that these conditions imply $\xi_1 \geq \xi_2$.)
The highest root is $\delta = \epsilon_1 - \epsilon_6$.

Write
$\lambda = n\epsilon + (\xi_1,\ldots,\xi_6)$,
where $\sum_i \xi_i = 0$. 
Assume that  $(\lambda,\delta)  = \xi_1 - \xi_6 \leq 2$.  
Combining the inequalities 
$\xi_1 - \xi_6  \leq 2$,
$\xi_4 \geq \xi_5$,  $\xi_4 \geq \xi_6$, and 
$\xi_1 + \xi_5 + \xi_6 \geq \frac{n}{2}$,
we see that  $\xi_4 \geq \frac{n-4}{6}$.
Since also $\xi_2 + \xi_3 + \xi_4 \leq 0$,  
$\xi_2 \geq \xi_4$, and $\xi_3 \geq \xi_4$,  we have 
$\xi_4 \leq 0$.
Moreover, in both cases, if the final inequality in the
sentence is an equality, so are all the preceding ones.
Since $n \geq 0$, $0 \geq \xi_4 \geq -\frac{4}{6}$.
Since $\lambda \in \tell$, $\xi_4 = 0$ or $\xi_4 = -\frac{1}{2}$.
In the former case, $\xi_2 = \xi_3 = \xi_4 = 0$, so
$\xi_1 + \xi_5 + \xi_6 = 0$, so $n = 0$.
Hence, $\lambda = (\epsilon_1 - \epsilon_6).$
In the latter case, $n$ is odd, so $\xi_4 = -\frac{1}{2} \geq  \frac{n-4}{6}$
implies that $n = 1$. 
In this case, $\lambda =  (\epsilon_1 - \epsilon_6)  
+ (\epsilon + \epsilon_1 + \epsilon_2  + \epsilon_6)$.
This proves claim (a).

Since $\delta = (\epsilon_1 -\epsilon_i) + (\epsilon_i - \epsilon_6),$
the roots $\pm (\epsilon_1 - \epsilon_i)$ and $\pm(\epsilon_i - 
\epsilon_6)$
lie in $L$ for all $1  < i < 6$.
Moreover, $\delta = (\epsilon + \epsilon_1 + \epsilon_i + \epsilon_j)
- (\epsilon + \epsilon_i + \epsilon_j +\epsilon_6)$,
so the roots $\pm(\epsilon + \epsilon_1 + \epsilon_i + \epsilon_j)$
and
$\pm(\epsilon + \epsilon_i + \epsilon_j + \epsilon_6)$
lie in $L$  for all $1 < i < j < 6$.
Since, for example, $\epsilon + \epsilon_1 + \epsilon_2 + \epsilon_3 =
\epsilon - \epsilon_4 - \epsilon_5 - \epsilon_6$, 
it follows easily that  $L$ contains all roots.

Let $\alpha \in \tell^*$ be a nonzero weight.
Write $\lambda = n\eps + \xi$ as before.
By applying the Weyl group, we may assume that
$\alpha = m\epsilon + (\zeta_1,\ldots,\zeta_6)$ is in the positive Weyl 
chamber.
Since $(\alpha,\lambda) \geq (\zeta,\xi)$, it is enough
to show that $(\zeta,\xi) > 1$.
This fact now follows
from the argument from $A_5$,
since $(\delta,\la)  = (\de,\xi)> 2$, 
since the
Weyl group contains the permutation group $S_5$,
and since $\zeta$ must 
satisfy $\zeta_i - \zeta_j \in \Z$. 

Finally,  $\delta$ is the only root in the positive Weyl chamber.

\MS

\NI{\bf (V, b)}\,\,
Recall that for the group $E_7$, $\ft = 
\ft^* = \{ \xi  \in \R^8 \mid \sum \xi_i = 0 \}$,
the roots are $\epsilon_i - \epsilon_j$
and $\epsilon_i + \epsilon_j + \epsilon_k + \epsilon_l$ for $i,j,k,$ and $l$
distinct, and the integral lattice is
$\tell = \{\xi \in \ft \mid \xi_i \pm \xi_j \in \Z \ \forall \ i,j \}.$
The positive Weyl chamber is 
$$
\{\xi \in \ft \mid \xi_2 \geq \cdots \geq \xi_8 \ \mbox{and}\ \xi_1 + \xi_6
+ \xi_7 + \xi_8 \geq 0 \}.
$$
(Note that this automatically implies that $\xi_1 \geq \xi_2$.)
The highest 
root is $\delta = \epsilon_1 - \epsilon_8$.

Assume  that  $(\delta,\lambda) = \lambda_1 - \lambda_8 \leq 2$.
Combining the inequalities $\lambda_1 - \lambda_8 \leq 2$,  
$\lambda_1 + \lambda_6 + \lambda_7 + \lambda_8 \geq 0$,
and $\lambda_5 \geq \lambda_i$ for $i = 6,7$ and $8$, we see
that $\lambda_5 \geq -\frac{1}{2}$. 
Since also $\lambda_2 + \lambda_3 + \lambda_4 + \lambda_5 \geq 0$
and $\lambda_i \geq \lambda_5$ for $i = 2,3$ and $4$,  $\lambda_5 \leq 0$. 
Moreover, in both cases, if the last inequality in the
sentence is an equality, all the inequalities are equalities.
Since $\lambda \in \ell$,
the only possibilities are $\lambda_5 = 0$ or $\lambda_5 = -\frac{1}{2}.$
In the former case, we must have $\lambda = \epsilon_1 - \epsilon_8$.
In the latter case, the only possibilities 
are
$\lambda = (\epsilon_1 + \epsilon_2 + \epsilon_3 + \epsilon_4) +
(\epsilon_1 - \epsilon_4)$, or
$\lambda = (\epsilon_1 + \epsilon_2 + \epsilon_3 + \epsilon_4) +
(\epsilon_1 - \epsilon_4) + (\epsilon_2 - \epsilon_3)$. 
The proves claim (a).

Since $\delta = (\epsilon_1 - \epsilon_i) + (\epsilon_i - \epsilon_8)$,
the roots $\pm(\epsilon_1 - \epsilon_i)$ and $\pm(\epsilon_i - \epsilon_8)$
lie in $L$ for all $1<i<8.$
Since $\delta = (\epsilon_1 + \epsilon_i + \epsilon_j + \epsilon_k)
- (\epsilon_i + \epsilon_j + \epsilon_k + \epsilon_8)$,
the roots
$\pm (\epsilon_1 + \epsilon_i + \epsilon_j + \epsilon_k)$
and
$\pm (\epsilon_i + \epsilon_j + \epsilon_k + \epsilon_8)$
also lie in $L$ for all $1<i<j<k<8$.
All roots can be written as a sum of these roots.  This proves claim (b).

Since $\tell \subset \Z^8 \cap \ft$, $\alpha_i - \alpha_j \in \Z$ for
every $\alpha \in \tell^* \subset \ft^*$.
Hence, the argument for claim (c) follows from the argument for $A_7$.

Finally,  $\delta$ is the only root in the positive Weyl chamber.

\MS

\NI{\bf (V, c)}\,\,
For the group $E_8$,
$\ft = \ft^* = \{\xi \in \R^9 \mid \sum \xi_i = 0$\}
and the roots are $\epsilon_i - \epsilon_j$, and
$\pm( \epsilon_i + \epsilon_j + \epsilon_k)$
for $i, j$ and  $k$ distinct.
Hence the integral lattice is
$$\tell = \{\xi \in \ft \mid 
3 \xi_i \in \Z \ \mbox{and}\  \xi_i - \xi_j \in \Z \ \forall \ i,j \}.$$
The positive Weyl chamber is 
$$
\{\xi \in \ft \mid \xi_2 \geq \cdots \geq \xi_9 \ \mbox{and} \ \xi_2
+ \xi_3 + \xi_4 \leq 0 \}.
$$
(Note that these conditions imply that $\xi_1 \geq \xi_2$.)
The highest root is $\delta = \epsilon_1 - \epsilon_9$.

Assume that  $(\delta,\lambda) = \lambda_1 - \lambda_9 \leq 2$. 
Combining the inequalities  
$$
\lambda_1 - \lambda_9 \leq 2,\quad
\lambda_1 + \lambda_5 + \lambda_6 + \lambda_7 + \lambda_8 + \lambda_9 \geq 0,
\quad
\lambda_i \leq \lambda_4,i > 4,
$$
we see that 
$\lambda_4 \geq - \frac{1}{3}$.
Since $\lambda_2 + \lambda_3 + \lambda_4 \leq 0$
and $\lambda_2 \geq \la_3\geq \lambda_4$,
$\lambda_4 \leq 0$.
Moreover, in both cases, if the last inequality in the
sentence is an equality, all the inequalities  
are equalities.
Since $\lambda \in \ell$, the only possibilities  are
$\lambda_4 = 0$ or $\lambda_4 = -\frac{1}{3}$.
In the former case, 
$\lambda = \epsilon_1 - \epsilon_9$.
In the latter case, 
$\lambda = (\epsilon_1 - \epsilon_9) + (\epsilon_1  + \epsilon_2 + \epsilon_9)$.
Claim (a) follows. 

We now notice that
$\delta = \epsilon_1 - \epsilon_9 = (\epsilon_1 - \epsilon_k)  +
(\epsilon_k - \epsilon_9) = (\epsilon_1 + \epsilon_i + \epsilon_j)
- (\epsilon_i + \epsilon_j + \epsilon_9)$ for all $1 < k < 9$ and $1 < i < j < 9$.
Therefore, the corresponding roots
$\pm(\epsilon_1 - \epsilon_k)$, $\pm(\epsilon_k - \epsilon_9)$,
$\pm(\epsilon_1 + \epsilon_i + \epsilon_j)$, and $\pm(\epsilon_i + \epsilon_j  +
\epsilon_9)$ all lie in $L$.
Since every root can be written as a sum of these roots, claim (b) follows. 

Since $\tell \subset \Z^9 \cap \ft$, $\alpha_i - \alpha_j \in \Z$ for
every $\alpha \in \tell^* \subset \ft^*$.
Hence, the argument for claim (c) carries over
from the argument for the group $A_8$.

Finally,  $\delta$ is the only root in the positive Weyl chamber.

\MS

\NI{\bf (VI)}\,\,
For the group $F_4$, $\ft = \ft^* = \R^4$.
The roots are $ \pm e_i$, $e_i \pm e_j$ for $i \neq j$, and
$\frac{1}{2}(\pm e_1 \pm e_2 \pm e_3 \pm e_4)$.
Hence the integral lattice is
$\tell = \{ \xi \in \Z^4 \mid \sum \xi_i \in 2 \Z \}$.
The positive Weyl chamber is
$$
\{\xi \in \ft \mid \xi_2 \geq \xi_3 \geq \xi_4 \geq 0 \ \mbox{and}
\ \xi_1 \geq \xi_2 + \xi_3 + \xi_4 \}.
$$
(Note that automatically $\xi_1 \geq \xi_2$.)
The highest root is $\delta = e_1 + e_2$.

The argument for claim (a) carries over
word for word from   the argument for $B_4$.

Notice that 
if $k=3$ or $4$  
\begin{eqnarray*}
\delta & = & e_1 + e_2 = (e_1) + (e_2)\; =\; (e_1 - e_k) + (e_2 + e_k)
\;=\; (e_1 + e_k) + (e_2 - e_k) \\ 
&=&
\frac{1}{2}(e_1 + e_2 + e_3 + e_4)  + \frac{1}{2}(e_1 + e_2 - e_3 - 
e_4)\\
& = & \frac{1}{2}(e_1 + e_2 - e_3 + e_4)  + \frac{1}{2}(e_1 + e_2 + e_3 - e_4).
\end{eqnarray*}
Hence, the corresponding roots all lie in $L$.
Since every root can be written as the sum of these roots, this
proves claim (b).

Since $-\id$ lies in the Weyl group,  we are done.

\MS

\NI{\bf (VII)}\,\,
For the group $G_2$,  $\ft = \ft^* = \{\xi \in \R^3 \mid \sum \xi_i = 0$\}.
The roots are $\pm \epsilon_i$ and $\epsilon_i - \epsilon_j$ 
for $i$ and  $j$ distinct.
The positive Weyl chamber is
$\{ \xi \in \ft^* \mid 0  \geq \xi_2 \geq \xi_3  \}$. 
(Note that automatically $\xi_1 \geq \xi_2$.)
The integral lattice  is $\tell =\Z^3 \cap \ft$.
The highest root is $\delta = \epsilon_1 - \epsilon_3$.

The argument for claim (a) follows 
the argument for $A_3$ word for word.

Since $\delta = (\epsilon_1 - \epsilon_2) + (\epsilon_2 - \epsilon_3)
= (\epsilon_1) + (-\epsilon_3)$,
the roots $\pm (\epsilon_1  - \epsilon_2)$, $\pm (\epsilon_2 - \epsilon_3)$,
$\pm \epsilon_1$ and $\pm \epsilon_3$ all lie in $L$.
Since every root can be written as a sum of these roots,  claim (b) follows.  

Since $-\id$ lies in the Weyl group, we are done.

\MS

\end{document}